\documentclass[11pt]{scrarticle}

\usepackage[T1]{fontenc}
\usepackage[utf8]{inputenc}
\usepackage{amsfonts,amsmath,amssymb,amsthm,mathtools,dsfont,
fullpage,xcolor,hyperref,indentfirst,microtype,authblk}
\usepackage[shortlabels]{enumitem}

\hypersetup{colorlinks,linkcolor={red!50!black},
citecolor={blue!50!black}}

\newcommand{\authorinfo}[3]
{\bigskip
\noindent\textsc{#1:}\\[5pt]#2\\[5pt]
\noindent\textit{Email address: }\href{mailto:#3}{#3}}

\newcommand{\eps}{\varepsilon}

\renewcommand{\P}{\mathds P}
\newcommand{\E}{\mathds E}
\newcommand{\R}{\mathds R}
\newcommand{\N}{\mathds N}
\renewcommand{\O}{\mathcal O}
\newcommand{\A}{\mathcal A}
\newcommand{\B}{\mathcal B}
\newcommand{\EE}{\mathcal E}
\newcommand{\T}{\mathcal T}
\newcommand{\hs}{\hspace{1pt}}
\newcommand{\e}{\mathrm e}
\newcommand{\dd}{\mathrm d}

\theoremstyle{plain}
\newtheorem{theorem}{Theorem}

\theoremstyle{definition}
\newtheorem{remark}{Remark}

\title{Explosion and non-explosion in pure birth Crump--Mode--Jagers branching processes}

\author{\textsc{Oleksii Galganov,
Andrii Ilienko}\thanks{Corresponding author. Supported by the Swiss National Science Foundation, grant no.~229505.}}

\date{}

\begin{document}
	
\maketitle

\begin{abstract}
In this short note, we provide an explicit sufficient condition for non-explosion of Crump--Mode--Jagers branching processes with pure birth reproduction. It shows that the standard sufficient condition for explosion, namely the convergence of the series of reciprocals of the birth rates, is --- at least for rate sequences without excessive oscillations --- remarkably close to being necessary. At the same time, it is not necessary in full generality: we construct a counterexample which also yields a general preferential attachment tree without fitness with an infinite path and no vertices of infinite degree, thereby answering an open question previously raised in the literature.
\end{abstract}


	
\section{Introduction}

Crump--Mode--Jagers (CMJ) branching processes, that is, branching processes constructed from general reproduction point processes (see, e.g., Chapter~6 in~\cite{J75}), form a very broad class. They are defined as follows. The initial ancestor is born at time $0$. Its children are born at times given by the atoms of a reproduction point process $\xi$. Each child then evolves in the same way, using an independent copy of $\xi$ shifted to start at the child’s birth time, and the construction continues recursively. In particular, writing $\delta_a$ for the unit point mass at $a$, if $\xi=N\delta_T$ for random $\N_0$-valued $N$ and $\R_+$-valued $T$, we obtain the Sevastyanov splitting process. Special cases include the Bellman--Harris process when $N$ and $T$ are independent, and the Galton--Watson process  when $T$ is constant. Many further examples can be found in~\cite{K16}. We do not describe the construction above in full formal detail, referring the reader instead to any of the recent papers \cite{I24,IL23,K16,L24} if needed.

In this note, we consider a pure birth CMJ branching process, that is, a CMJ process whose reproduction $\xi$ is given by a pure birth point process with rates $\lambda_i$, $i\ge1$.
Recall that a pure birth counting process is a continuous-time Markov chain $(X_t)_{t\ge0}$ on $\N$, started from $X_0=1$\footnote{The state space is usually taken to be $\N_0$ and the process is started at $0$. We adopt a shift by one for notational convenience.},
with infinitesimal generator matrix $G=(g_{ij})_{i,j\in\N}$ given by
\[
g_{i,i+1}=\lambda_i,\qquad g_{ii}=-\lambda_i,\qquad g_{ij}=0\text{ for }j\notin\{i,i+1\},\quad i\in\N;
\]
see, e.g., Section~2.5 in~\cite{N98}. The associated pure birth point process is $\xi=\sum_{i=1}^\infty\delta_{T_i}$,
where $T_i=\inf\{t\ge0:\,X_t>i\}$. It is well known that
\[(T_i)_{i\ge1}\overset d=\Biggl(\sum_{j=1}^iE_j\Biggr)_{i\ge1},\]
where $E_j$, $j\ge1$, are independent $\mathsf{Exp}(\lambda_j)$ random variables.

Our interest in pure birth CMJ processes stems from the fact that they play a central role in the study of preferential attachment random recursive trees. We recall the basic definitions. Let $f:\N_0\to(0,\infty)$ be a fixed attachment function. We now construct a consistent sequence $(\T_n)_{n\in\N_0}$ of random trees as follows. The tree $\T_0$ consists of a single vertex $0$. At step $n\ge1$, a new vertex $n$ appears and an edge is added from an existing vertex $m$ to $n$, where $m$ is chosen with probability proportional to $f$ of its current out-degree:
\begin{equation}
\label{eq:tree}
\P\{m\to n\mid\T_{n-1}\}=\frac{f(\deg_{n-1}^+(m))}{\sum_{j=0}^{n-1}f(\deg_{n-1}^+(j))},\qquad m=0,\ldots,n-1.
\end{equation}
Here the out-degree $\deg_{n-1}^+(m)$ is the number of children of $m$ in $\T_{n-1}$, that is, immediately before $n$ is attached. One also considers preferential attachment trees with fitness, where $f$ in \eqref{eq:tree} depends not only on the out-degree of $m$ but also on its random weight $W_m$, with $(W_m)_{m\in\N_0}$ forming an i.i.d.\ sequence; see, e.g., \cite{I24,IL23,L24}.

The connection between preferential attachment trees and CMJ processes was first observed in \cite{P94} and was used there to study the asymptotics of the height of trees with affine preferential attachment. Subsequently, this approach has become standard for affine and general preferential attachment trees; see, e.g., \cite{IK18,OS05,RTV07}. The core of this connection, going back to the classical poissonization technique, is as follows. Let $\tau_n$, $n\ge0$, be the birth time of the $n$th individual in the pure birth CMJ process with rates $\lambda_i=f(i-1)$, $i\in\N$. Then the sequence of genealogical trees of this process, observed at times $\tau_n$, after relabelling the vertices in each tree in chronological order of birth, has the same distribution as $(\T_n)_{n\in\N_0}$. This remains true for preferential attachment trees with fitness as well (see Proposition~3.14 in \cite{I24}), but in this case one has to replace the pure birth CMJ process by a CMJ process that is pure birth conditionally on the random weights.

Recall that explosion in a CMJ process is the event $\EE=\{\lim_{n\to\infty}\tau_n<\infty\}$, that is, the birth of infinitely many individuals in finite time. It follows from Proposition 3.5 in \cite{I24} that for a pure birth CMJ process one has $\P(\EE)\in\{0,1\}$. In view of the above discussion, the question of necessary and sufficient conditions for $\P(\EE)=1$ becomes particularly important. Specifically, the presence or absence of explosion in the associated pure birth CMJ process determines the shape of the corresponding limit preferential attachment tree without fitness --- namely, whether it contains infinite stars, that is, vertices of infinite degree, and/or infinite paths. Unfortunately, deriving necessary and sufficient conditions for explosion in a reasonably explicit and usable form appears to be beyond current methods. This stands in sharp contrast to, say, Bellman--Harris branching processes, for which such conditions are known in many cases; see, e.g., \cite{G73}. The general operator approach, presented in its most complete form in Section 2 of \cite{K16}, does not yield any explicit conditions either.

In this situation, one can only rely on sufficient conditions for explosion. A simple condition stated directly in terms of the rates, namely
\begin{equation}
\label{eq:conv}
\sum_{i=1}^\infty \lambda_i^{-1}<\infty,
\end{equation}
follows immediately from the standard explosion criterion for the non-branching pure birth process. In part (ii) of Theorem \ref{th:main}, we give a fairly general sufficient condition for non-explosion. Remark \ref{rem:log} shows that, for rate sequences that do not oscillate too strongly, the non-explosion criterion in part (ii) nearly complements \eqref{eq:conv}; in this sense, \eqref{eq:conv} is close to being necessary for explosion. Part (iii), however, shows by means of a rather pathological sequence of rates that it is not exactly necessary. This yields an example of a general preferential attachment tree without fitness with no infinite stars and a unique infinite path. In turn, this answers the open question raised in Section~3.4 of~\cite{IL23} concerning whether such a situation can occur.
	
\section{Main result and discussion}
	
The following theorem provides sufficient conditions for explosion/non-explosion of a pure birth CMJ process. For completeness, we include \eqref{eq:conv} as (i).

\begin{theorem}~
\label{th:main}
\begin{enumerate}[\normalfont(i)]
	\item If $\sum_{i=1}^\infty \lambda_i^{-1}<\infty$, then the process is explosive.
	\item If
	\begin{equation}
	\label{eq:lim1}
	\lim_{i\to\infty}\frac{\lambda_i}i=\infty
	\end{equation}
	and
	\begin{equation}
	\label{eq:lim2}
	\liminf_{n\to\infty}\!\!\sum_{i=\lceil\eps\log n\rceil}^n\!\! \lambda_i^{-1}>0\quad\text{for some $\eps>0$},
	\end{equation}
	then the process is non-explosive.
	
	The same holds if $(\lambda_i)$ is dominated by a sequence satisfying \eqref{eq:lim1} and \eqref{eq:lim2}.
	\item There exists a sequence $(\lambda_i)$ such that $\sum_{i=1}^\infty\lambda_i^{-1}=\infty$ and the process is explosive.
\end{enumerate}
\end{theorem}

Condition \eqref{eq:lim2} may be viewed in light of the following elementary observation. It is easy to see that
$\sum_{i=1}^\infty\lambda_i^{-1}=\infty$ alone implies
$\liminf_{n\to\infty}\sum_{i=b_n}^n\lambda_i^{-1}>0$
for some integer-valued sequence $b_n\to\infty$ growing sufficiently slowly. Condition \eqref{eq:lim2} requires that such a lower cutoff $b_n$ can be chosen to be logarithmic in $n$. The following remark shows that the class of such rate sequences is very broad.

\begin{remark}
\label{rem:log}
Consider the class of sequences $(\lambda_i)$ that may come very close to the convergence boundary of the series $\sum_{i=1}^\infty \lambda_i^{-1}$:
\[\lambda_i=\O\bigl(i\cdot\log i\cdot\log\log i\cdot\ldots\cdot\log^{(k)} i\bigr)\quad\text{for some }k\ge1,\]
where $\log^{(k)}$ denotes the $k$-fold iterated logarithm. Since, as $c>0$ and $n\to\infty$,
\begin{equation}
\label{eq:log}
\sum_{i=\lceil\log n\rceil}^n\Bigl(c\hs i\cdot\log i\cdot\log\log i\cdot\ldots\cdot\log^{(k)} i\Bigr)^{-1}\!\!\sim c^{-1}\bigl(\log^{(k+1)}n-\log^{(k+2)}n\bigr)\to\infty,
\end{equation}
the process is non-explosive by (ii). The substantial gap between \eqref{eq:lim2} and \eqref{eq:log} leaves ample room for constructing non-explosive sequences of rates much closer to the convergence boundary.

This example shows that the sufficient condition for explosion in (i) is, in a sense, very close to being necessary, whereas (iii) demonstrates that it is not exactly so.
\end{remark}

\begin{remark}
\label{rem:Iyer}
	After this note had already been submitted, Tejas Iyer kindly pointed out to us that a further strong sufficient condition for non-explosion can be derived from the results of his recent preprint \cite{I24p} on persistent hubs in CMJ processes. A persistent hub is a vertex which, from some time onward, remains a vertex of maximal degree.
	
	More precisely, consider a pure birth CMJ process whose birth rates satisfy
	\begin{gather}
		\label{eq:ph1}
		\sum_{i=1}^\infty \lambda_i^{-1}=\infty,\\
		\max_{i\le n}\frac{\lambda_i}{i}\le\kappa\frac{\lambda_n}{n}
		\label{eq:ph2}
	\end{gather}
	for all $n$ and some constant $\kappa>0$. Under condition \eqref{eq:ph2}, Theorem 2.16 of \cite{I24p} implies that the process has a unique persistent hub a.s. (Note that the function $f(i)$ in \cite{I24p} corresponds to our $\lambda_{i+1}$.) If the process were explosive, then, by Lemma 3.6 of \cite{I24p}, the degree of the persistent hub at the explosion time would be infinite. This is impossible in view of \eqref{eq:ph1}. Thus, \eqref{eq:ph1} and \eqref{eq:ph2} yield another sufficient condition for non-explosion.
	
	The conditions \eqref{eq:lim1}+\eqref{eq:lim2} and \eqref{eq:ph1}+\eqref{eq:ph2} are not comparable: one can easily construct examples in which the former applies but the latter does not, and conversely. The former allows for moderate oscillations in the birth rates, at the cost of the additional condition \eqref{eq:lim2}. The latter is restricted to rather regular rates, in the sense of \eqref{eq:ph2}, but, within this class, requires nothing beyond \eqref{eq:ph1}.
	
	Finally, let us note that, by strengthening condition \eqref{eq:lim1}, for instance to
	\[\lim_{i\to\infty}\frac{\lambda_i}{i\log i}=\infty
	\qquad\text{or}\qquad
	\lim_{i\to\infty}\frac{\lambda_i}{i\log i\log\log i}=\infty,\]
	one can use a similar reasoning to lower the starting index in the summation in \eqref{eq:lim2}, thereby partially bridging the gap between \eqref{eq:lim1}+\eqref{eq:lim2} and \eqref{eq:ph1}+\eqref{eq:ph2}. However, a general formulation of this kind would be considerably less transparent.	
\end{remark}

\begin{remark}
\label{rem:tree}
Consider the preferential attachment tree without fitness, which is the discrete-time skeleton of the explosive CMJ process in (iii). By item 3 of Theorem 3.16 in \cite{I24}, this tree has a unique infinite path and no infinite stars, that is, vertices of infinite degree. A related example was constructed in Corollary~3.19 of~\cite{I24} in the context of linear preferential attachment trees with heavy-tailed fitness, in particular in the Bianconi--Barab\'asi model. (iii) shows that, surprisingly, the same effect can occur for some general preferential attachment trees with no fitness at all. As already mentioned, this answers the open question raised in Section 3.4 of \cite{IL23}.

Thus, all three possible shapes of preferential attachment trees that are known to arise in the presence of fitness --- (a) every vertex is an infinite star and every path is infinite, (b) a unique infinite star and no infinite paths, and (c) a unique infinite path and no infinite stars --- can also be obtained without fitness: (a) occurs whenever the associated CMJ process is non-explosive; (b) is realized, for instance, in super-linear power preferential attachment (Theorem~1.2 in~\cite{OS05}) as well as under the assumptions of Corollary~2.16 in~\cite{L24}; and (c) corresponds to the explosive CMJ process in (iii).
\end{remark}

\begin{remark}
The rates $\lambda_i$ used in the proof of (iii) are highly irregular: they take astronomically large values along long stretches, while returning to the baseline value $1$ infinitely often, at an increasingly sparse set of indices. In particular, they satisfy neither \eqref{eq:lim1} nor \eqref{eq:lim2}; \eqref{eq:ph2} fails as well. We have not been able to construct an example with, say, monotone or regularly varying rates. It would be interesting to determine whether such a construction is possible at all.
\end{remark}

\subsection*{Applications to preferential attachment trees}
Consider a consistent sequence $(\T_n)_{n\in\N_0}$ of general preferential attachment trees without fitness, with attachment probabilities given by \eqref{eq:tree}, and the corresponding limit tree $\bigcup_{n\in\N_0}\!\T_n$. By a consequence of K\H onig's infinity lemma (see, e.g., Proposition 8.2.1 in \cite{D25}), any infinite tree contains either an infinite star, that is, a vertex with infinitely many children, or an infinite path, that is, a ray $0\to v_1\to v_2\to\ldots$ For preferential attachment random trees, possibly with fitness, the picture is even sharper. If the corresponding CMJ process explodes, then, by Theorem 2.12 in \cite{IL23}, the limit tree has either a unique infinite star or a unique infinite path, but not both. If, on the other hand, the CMJ process does not explode, then each vertex keeps producing children indefinitely; thus, every vertex of the limit tree is an infinite star, and every finite path can be extended to an infinite one.	
	
Consider now a preferential attachment tree without fitness, which is the discrete-time skeleton of the explosive CMJ process in (iii). By item 3 of Theorem 3.16 in \cite{I24}, this tree has a unique infinite path and no infinite stars. A related example was constructed in Corollary~3.19 of~\cite{I24} in the context of linear preferential attachment trees with heavy-tailed fitness, in particular in the Bianconi--Barab\'asi model. (iii) shows that the same effect can occur for some general preferential attachment trees with no fitness at all. As already mentioned, this answers the open question raised in Section 3.4 of \cite{IL23}.
	
Thus, all three possible shapes of limit preferential attachment trees that are known to arise in the presence of fitness can also occur without fitness:
\begin{enumerate}[(a)]
	\item every vertex is an infinite star, and every finite path can be extended to an infinite one --- whenever the associated CMJ process is non-explosive;
	\item a unique infinite star and no infinite paths --- for instance, for super-linear power preferential attachment (see Theorem~1.2 in \cite{OS05}) as well as under the assumptions of Corollary~2.16 in \cite{L24};
	\item a unique infinite path and no infinite stars --- in the tree corresponding to the explosive CMJ process in (iii).
\end{enumerate}	
	
Finally, we note that a modification of the example in (iii) was very recently used by Tejas Iyer to disprove the conjecture that the condition
$\sum_{j=0}^\infty\frac1{f(j)^2}<\infty$
is sufficient for the existence of a persistent hub in preferential attachment trees without fitness \cite{I26}.

\section{Proof of Theorem \ref{th:main}}
As noted above, (i) follows immediately from the standard explosion criterion for a usual pure birth process, since in this case even the initial individual produces infinitely many offspring in finite time.

For (ii), first fix $\eps>0$. Since
\[\lim_{\beta\to\infty}\lim_{\substack{t\to0,\\r\to\infty}}\Bigl(t+\frac 2\beta\Bigr)
\Bigl(1+\frac{\beta}{r\eps}\Bigr)=0,\]
by \eqref{eq:lim2} we may choose $t,\beta,r>0$ such that
\begin{equation}
\label{eq:in1}
\liminf_{n\to\infty}\!\!\sum_{i=\lceil\eps\log n\rceil}^n\!\!\lambda_i^{-1}>
\Bigl(t+\frac 2\beta\Bigr)\Bigl(1+\frac{\beta}{r\eps}\Bigr).
\end{equation}
Moreover, by \eqref{eq:lim1}, $i\ge\lceil\eps\log n\rceil$ implies $\lambda_i\ge ri\ge r\eps\log n$ for large $n$, and hence $\beta\log n\le\frac{\beta \lambda_i}{r\eps}$. Therefore, by \eqref{eq:in1},
\begin{equation}
\label{eq:in2}
\sum_{i=1}^n\frac 1{\lambda_i+\beta\log n}\ge\sum_{i=\lceil\eps\log n\rceil}^n\frac 1{\lambda_i+\beta\log n}\ge\Bigl(1+\frac{\beta}{r\eps}\Bigr)^{-1}
\!\!\!\sum_{i=\lceil\eps\log n\rceil}^n\!\lambda_i^{-1}>t+\frac 2\beta
\end{equation}
for large $n$.

It is well known that, to prove non-explosion of a CMJ process constructed from a reproduction point process $\xi$ with no atom at the origin, it suffices to show that $\xi$ has a finite intensity measure in a neighborhood of the origin, that is,
\begin{equation}
\label{eq:mean}
\E\xi\bigl([0,t]\bigr)<\infty
\end{equation}
for some $t>0$; see Theorem 6.2.2 in \cite{J75} or Theorem 3.1\hs(b) in \cite{K16}.

Let $E_i\sim\mathsf{Exp}(\lambda_i)$, $i\ge1$, be independent. Then, in our case,
\begin{equation}
\label{eq:series}
\begin{aligned}
\E\xi\bigl([0,t]\bigr)=\sum_{n=1}^\infty\P\biggl\{\sum_{i=1}^nE_i\le t\biggr\}=\sum_{n=1}^\infty\P\biggl\{\prod_{i=1}^n\e^{-\beta E_i\log n}\ge\e^{-\beta t\log n}\biggr\}\\\le
\sum_{n=1}^\infty\e^{\beta t\log n}\prod_{i=1}^n\E\e^{-\beta E_i\log n}=
\sum_{n=1}^\infty\e^{\beta t\log n}\prod_{i=1}^n
\frac{\lambda_i}{\lambda_i+\beta\log n}.
\end{aligned}
\end{equation}
The $n$th term on the right-hand side is
\[\exp\biggl\{\beta t\log n+\sum_{i=1}^n\log\Bigl(1-\frac{\beta\log n}{\lambda_i+\beta\log n}\Bigr)\biggr\}\le\exp\Bigl\{\beta\log n\Bigl(t-\sum_{i=1}^n\frac{1}{\lambda_i+\beta\log n}\Bigr)\biggr\},\]
which, in view of \eqref{eq:in2}, is bounded by $n^{-2}$ for large $n$.
Hence, the series on the right-hand side of \eqref{eq:series} converges, and \eqref{eq:mean} follows.

If $(\lambda_i)$ is only dominated by a sequence satisfying \eqref{eq:lim1} and \eqref{eq:lim2}, it suffices to increase all $E_i$ accordingly, which does not cause explosion. This completes the proof of (ii).

We now turn to the construction of an example for~(iii). By Theorem~3.8 in~\cite{I24}, to prove explosion of a CMJ process with reproduction point process $\xi$, it suffices to show that
\[\sum_{i=1}^\infty\bigl(1-\P\bigl\{\xi\bigl([0,t_i]\bigr)>M_{i+1}\bigr\}\bigr)^{M_i}<\infty\]
for some $t_i>0$ with $\sum_{i=1}^\infty t_i<\infty$ and some $M_i\in\N$. Since $\sum_{i=1}^\infty\e^{-i}<\infty$ and $\log(1-x)\le-x$, it is enough to prove that
\begin{equation}
\label{eq:main}
\lim_{i\to\infty}\frac{M_i}i\cdot\P\bigl\{\xi\bigl([0,t_i]\bigr)>M_{i+1}
\bigr\}=\infty.
\end{equation}

Consider a pure birth CMJ process with rates $\lambda_i$ given by
\begin{equation}
\label{eq:rates}
1,\underbrace{\alpha_1,\ldots,\alpha_1}_{d_1\text{ times}},
1,\underbrace{\alpha_2,\ldots,\alpha_2}_{d_2\text{ times}},1,\ldots
\end{equation}
with $d_k=4^{k^2}$ and $\alpha_k=2^{2^{k^3}}$. Clearly, $\sum_{i=1}^\infty\lambda_i^{-1}=\infty$. Set $t_i=2^{-i}$ and $M_i=i+\sum_{k=1}^i d_k$, so that $M_i>4^{i^2}$. For independent $E_j\sim\mathsf{Exp}(\lambda_j)$, we have:
\begin{multline*}
\P\bigl\{\xi\bigl([0,t_i]\bigr)>M_{i+1}\bigr\}=\P\Biggl\{\sum_{j=1}^{M_{i+1}+1}E_j\le t_i\Biggr\}\\=\P\Biggl\{E_1+\!\sum_{j=2}^{1+d_1}E_j+E_{2+d_1}+
\!\!\sum_{j=3+d_1}^{2+d_1+d_2}\!E_j+\ldots+
\!\!\!\sum_{j=i+2+\sum_{k=1}^id_k}^{i+1+\sum_{k=1}^{i+1}d_k}\!\!E_j+
E_{i+2+\sum_{k=1}^{i+1}d_k}\le 2^{-i}\Biggr\}.
\end{multline*}
In view of~\eqref{eq:rates}, all single summands above are $\mathsf{Exp}(1)$, the terms in the first sum are $\mathsf{Exp}(\alpha_1)$, and so on, with the terms in the last sum being $\mathsf{Exp}(\alpha_{i+1})$. Denote the sum of all single summands by $S_0$, and the $k$th sum by $S_k$, $k=1,\ldots,i+1$. Then $S_0\sim\mathsf{\Gamma}(i+2,1)$, $S_k\sim\mathsf{\Gamma}(d_k,\alpha_k)$, and all these are independent. Hence,
\begin{equation}
\label{eq:S}
\P\bigl\{\xi\bigl([0,t_i]\bigr)>M_{i+1}\bigr\}\ge\P\bigl\{S_0\le 2^{-i-1}\bigr\}\prod_{k=1}^{i+1}\P\bigl\{S_k\le2^{-k-i-1}\bigr\}.
\end{equation}
Let $X\sim\mathsf{\Gamma}(d,\alpha)$ with $d\in\N$. Then, for any $t>0$, writing $\alpha t=y$, we have:
\begin{equation}
\label{eq:Gamma}
\P\{X\le t\}
=\frac{\alpha^d}{(d-1)!}\int_0^t x^{d-1}\e^{-\alpha x}\,\dd x
\ge\frac{\alpha^d\e^{-y}}{(d-1)!}\int_0^t x^{d-1}\,\dd x
=\frac{y^d\e^{-y}}{d\hs!}\ge\e^{-y}\Bigl(\frac yd\Bigr)^d.
\end{equation}
Therefore,
\[\P\bigl\{S_0\le 2^{-i-1}\bigr\}\ge\e^{-2^{-i-1}}
\biggl(\frac{2^{-i-1}}{i+2}\biggr)^{i+2},\]
and a simple calculation shows that the right-hand side dominates
$i\hs 3^{-i^2}$ as $i\to\infty$.
The remainder of the proof is devoted to showing that
\begin{equation}
\label{eq:mainprod}
\lim_{i\to\infty}\biggl(\frac 43\biggr)^{\!i^2}\!\!\cdot\prod_{k=1}^{i+1}\P\bigl\{S_k\le2^{-k-i-1}\bigr\}
=\infty,
\end{equation}
which, combined with the previous bound and the inequalities~\eqref{eq:S} and $M_i>4^{i^2}$, yields~\eqref{eq:main}.

Let $t_{i,k}=2^{-k-i-1}$ and $y_{i,k}=\alpha_kt_{i,k}=2^{2^{k^3}-k-i-1}$.
For $k=1$, we have $d_1=\alpha_1=4$ and $y_{i,1}=2^{-i}\le1$, and hence, by \eqref{eq:Gamma},
\begin{equation}
\label{eq:p1}
\P\{S_1\le t_{i,1}\}\ge\e^{-1}2^{-4i-8}.
\end{equation}
The remaining factors will be bounded in two different ways: for large $k$ by Markov's inequality, and for small $k$ by~\eqref{eq:Gamma}.
Define
\begin{equation}
\begin{gathered}
\label{eq:B}
\A_i=\bigl\{k\in\{2,\ldots,i+1\}: 2^{k^3}\ge 2k^2+k+i+2\bigr\},\\
\B_i=\bigl\{k\in\{2,\ldots,i+1\}: 2^{k^3}\le 2k^2+k+i+1\bigr\}.
\end{gathered}
\end{equation}
For $k\in\A_i$, we have
$\E S_k=\frac{d_k}{\alpha_k}=2^{2k^2-2^{k^3}}$, and hence, by Markov's inequality and the definition of~$\A_i$,
\begin{equation}
\label{eq:p2}
\P\{S_k\le t_{i,k}\}\ge1-\frac{\E S_k}{t_{i,k}}
=1-2^{-2^{k^3}+2k^2+k+i+1}\ge\frac12.
\end{equation}
Since $2^{k^3}-2k^2-k-1\ge0$ for all $k\ge2$, the definition of $\mathcal B_i$ gives
\[-i\le 2^{k^3}-2k^2-k-i-1\le0,\qquad k\in\mathcal B_i,\]
which implies $2^{-i}\le\frac{y_{i,k}}{d_k}\le1$.
Thus, by \eqref{eq:Gamma},
\begin{equation}
\label{eq:p3}
\P\{S_k\le t_{i,k}\}\ge\e^{-d_k}\,2^{-i d_k}.
\end{equation}
Let $\lvert\cdot\rvert$ denote the cardinality of a set and $D_i=\sum_{k\in\B_i} d_k$. Combining \eqref{eq:p1}, \eqref{eq:p2}, \eqref{eq:p3}, we obtain
\begin{equation}
\label{eq:f1}
\prod_{k=1}^{i+1}\P\bigl\{S_k\le2^{-k-i-1}\bigr\}\ge 2^{-|\A_i|-iD_i-4i-8}\e^{-D_i-1}
\ge2^{-iD_i-5i-8}\e^{-D_i-1}.
\end{equation}
From~\eqref{eq:B} it is apparent that $\lvert\B_i\rvert$ is of order $(\log_2 i)^{1/3}$. More rigorously, elementary but routine calculations show that
$\lfloor(\log_2 i)^{1/3}\rfloor+2\notin\B_i$ for large $i$; thus, $\lvert\B_i\rvert<(\log_2 i)^{1/3}$ for such $i$. Hence, since $d_k$ are increasing, we have
\begin{equation}
\label{eq:f2}
D_i\le d_{\lfloor(\log_2 i)^{1/3}\rfloor+1}\hs\lvert\B_i\rvert<4^{\left((\log_2 i)^{1/3}+1\right)^2}
(\log_2 i)^{1/3}.
\end{equation}
A straightforward verification shows that the right-hand side of~\eqref{eq:f1}, together with the bound~\eqref{eq:f2}, satisfies~\eqref{eq:mainprod}. This completes the proof of (iii).
\bigskip

\noindent{\large\bfseries Acknowledgments.} The authors are grateful to Alexander Iksanov and J\'ulia Komj\'athy for helpful and stimulating discussions, and to Tejas Iyer for bringing to our attention another sufficient condition for non-explosion, given in Remark \ref{rem:Iyer}.

\bigskip

\authorinfo{Oleksii Galganov}
{Department of Mathematical Analysis and Probability Theory,\\Igor Sikorsky Kyiv Polytechnic Institute, Prospect Beresteiskyi~37, 03056, Kyiv, Ukraine}{galganov.oleksii@lll.kpi.ua}
\bigskip

\authorinfo{Andrii Ilienko}
{Department of Mathematical Analysis and Probability Theory,\\Igor Sikorsky Kyiv Polytechnic Institute, Prospect Beresteiskyi~37, 03056, Kyiv, Ukraine;\\[5pt]
Institute of Mathematical Statistics and Actuarial Science,\\University of Bern, Alpeneggstrasse~22, 3012, Bern, Switzerland.} {andrii.ilienko@unibe.ch}

\end{document}